\documentclass[11pt]{amsart}

\usepackage[T1]{fontenc}
\usepackage[utf8]{inputenc}
\usepackage[english]{babel}
\usepackage{microtype}

\usepackage{geometry}
\usepackage{mathtools}
\usepackage{amssymb}
\usepackage{mathrsfs}
\usepackage{bm}

\usepackage[table,dvipsnames]{xcolor}

\usepackage{graphicx}
\usepackage{tikz}
\usepackage{tikz-cd}
\usepackage{float}

\usepackage[colorlinks=true,linkcolor=blue,citecolor=blue,urlcolor=blue]{hyperref}

\usepackage{enumitem}
\setlist{leftmargin=*, itemsep=2pt, topsep=3pt}

\numberwithin{equation}{section}

\theoremstyle{plain}
\newtheorem{thm}{Theorem}[section]
\newtheorem{lemma}[thm]{Lemma}
\newtheorem{prop}[thm]{Proposition}
\newtheorem{cor}[thm]{Corollary}
\newtheorem{claim}[thm]{Claim}
\newtheorem{conjecture}[thm]{Conjecture}

\theoremstyle{definition}
\newtheorem{defn}[thm]{Definition}
\newtheorem{example}[thm]{Example}
\newtheorem{remark}[thm]{Remark}

\newtheorem{mainthm}{Theorem}

\newcommand{\btheorem}{\begin{thm}}
\newcommand{\etheorem}{\end{thm}}
\newcommand{\blemma}{\begin{lemma}}
\newcommand{\elemma}{\end{lemma}}
\newcommand{\bproposition}{\begin{prop}}
\newcommand{\eproposition}{\end{prop}}
\newcommand{\bcorollary}{\begin{cor}}
\newcommand{\ecorollary}{\end{cor}}
\newcommand{\bdefinition}{\begin{defn}}
\newcommand{\edefinition}{\end{defn}}
\newcommand{\bclaim}{\begin{claim}}
\newcommand{\eclaim}{\end{claim}}
\newcommand{\bproof}{\begin{proof}}
\newcommand{\eproof}{\end{proof}}
\newcommand{\bremark}{\begin{remark}}
\newcommand{\eremark}{\end{remark}}
\newcommand{\bexample}{\begin{example}}
\newcommand{\eexample}{\end{example}}

\newcommand{\PP}{\mathbb{P}}

\newcommand{\OO}{\mathcal{O}}

\title{Positivity of Exterior Powers of Tangent Bundles and Their Subsheaves}
\author{Yuting Liu}
\date{}

\begin{document}

\begin{abstract}
We study positivity properties of exterior powers of tangent bundles and of their subsheaves. We prove a structure theorem for projective manifolds with nef third exterior power of the tangent bundle. We also prove a rigidity theorem for ample subsheaves of $\wedge^2T_X$ under a rank condition. Finally we study regular foliations whose exterior powers are strictly nef.
\end{abstract}

\maketitle

\tableofcontents
\section{Introduction}

Positivity properties of tangent bundles play a central role in the classification theory of projective manifolds. Since Mori's proof of Hartshorne's conjecture \cite{mori1979projective} and the solution of the Frankel conjecture by Siu and Yau \cite{siu1980compact}, it has been understood that positivity assumptions on $T_X$ impose strong restrictions on the geometry of $X$. This philosophy has been developed in many directions. For instance, Mok \cite{mok1988uniformization} classified compact K\"ahler manifolds with semipositive holomorphic bisectional curvature, while Demailly, Peternell and Schneider \cite{demailly1994compact} established a structure theorem for compact K\"ahler manifolds with nef tangent bundle.

A natural refinement is to study positivity of exterior powers of the tangent bundle. The case of $\wedge^2T_X$ has been studied extensively. Campana and Peternell \cite{campana1991projective} classified projective threefolds with nef $\wedge^2T_X$, and Yasutake \cite{yasutake2012second} obtained related results for certain Fano fourfolds. Under stronger positivity assumptions, Cho and Sato \cite{cho1995smooth} proved that a projective manifold with ample $\wedge^2T_X$ is either projective space or a smooth hyperquadric. This was later generalized to the strictly nef setting by Li, Ou and Yang \cite{li2019projective}. More recently, Watanabe \cite{watanabe2021positivity} proved that if $\wedge^2T_X$ is nef and $\dim X\ge 3$, then either $T_X$ is nef or $X$ is Fano.

Our first main result concerns the next case, namely the nefness of $\wedge^3T_X$.

\begin{mainthm}\label{wedge3}
Let $X$ be a smooth projective manifold such that $\wedge^3T_X$ is nef and $\dim X\ge 4$. Then at least one of the following holds.

\noindent\textup{(1)} The tangent bundle $T_X$ is nef.

\noindent\textup{(2)} The manifold $X$ is Fano.

\noindent\textup{(3)} Up to a finite \'etale cover, $X$ admits a locally trivial Fano fibration over an elliptic curve.
\end{mainthm}

Theorem \ref{wedge3} shows that nefness of the third exterior power still imposes a strong global structure. Compared with Watanabe's theorem for $\wedge^2T_X$, one new phenomenon appears: a locally trivial Fano fibration over an elliptic curve may occur after a finite \'etale cover.

We then turn to positivity of subsheaves of exterior powers of tangent bundles. Andreatta and Wi\'sniewski \cite{andreatta2001manifolds} proved that if $T_X$ contains an ample locally free subsheaf, then $X$ is projective space. Araujo, Druel and Kov\'acs \cite{araujo2007cohomological} proved a cohomological characterization of projective spaces and hyperquadrics using sections of $\wedge^pT_X\otimes L^{-p}$. Motivated by these results, Kov\'acs proposed the following conjecture, recorded by Ross in \cite{ross2010characterizations}.

\begin{conjecture}[\cite{ross2010characterizations}]\label{conjecture}
Let $X$ be a projective manifold of dimension $n$, and let $\mathcal F$ be an ample vector bundle on $X$ which is a subsheaf of $\wedge^pT_X$ for some $p>0$. If $\mathcal F=\wedge^p\mathcal E$ for some ample vector bundle $\mathcal E$, then either $X\simeq \PP^n$, or $p=n$ and $X\simeq Q_n$.
\end{conjecture}

The rank-one case was proved by Druel and Paris \cite{druel2013characterizations}. Ross \cite{ross2010characterizations} proved the conjecture in dimension two and also when the Picard number is one. We prove the conjecture in the case $p=2$ under a natural rank condition.

\begin{mainthm}\label{maintheorem}
Let $X$ be a projective manifold with $\dim X=n\ge 3$, and let $\mathcal F$ be an ample vector bundle which is a subsheaf of $\wedge^2T_X$. Assume that $\mathcal F=\wedge^2\mathcal E$ for some ample vector bundle $\mathcal E$. If $\operatorname{rank}(\mathcal F)\ge n$, then $X\simeq \PP^n$.
\end{mainthm}

Finally, we study a weaker positivity condition for foliations. Strict nefness is weaker than ampleness, and therefore one cannot expect the same rigidity statements. Nevertheless, it still imposes strong restrictions on the geometry of the foliation. Liu, Ou and Yang \cite{liu2020projective} studied projective manifolds whose tangent bundle contains a strictly nef subsheaf. For foliations, Ou \cite{ou2021foliations} proved a structure theorem under the nefness of the anti-canonical class of a foliation. We prove the following statement for strictly nef exterior powers of regular foliations.

\begin{mainthm}\label{main}
Let $\mathcal F$ be a regular foliation of rank $d$ on a projective manifold $X$. If $\wedge^r\mathcal F$ is strictly nef for some $1\le r<d$, then there exists a locally constant fibration
\[
f:X\to Y
\]
such that the fibers are Fano manifolds and the base $Y$ is a canonically polarized hyperbolic projective manifold. Moreover, the foliation $\mathcal F$ is induced by the fibration $f$. In particular, $\mathcal F$ is algebraically integrable.
\end{mainthm}

The limiting case $r=d$ is precisely the condition that
$-K_{\mathcal F}=\det\mathcal F$ be strictly nef. For regular foliations,
this determinant case was treated in \cite[Corollary C]{liu2021algebraic}:
$\mathcal F$ is induced by a locally constant fibration with rationally
connected fibers over a canonically polarized hyperbolic projective manifold.
Theorem \ref{main} treats the complementary range $1\le r<d$, where the
fibers are shown to be Fano.

As a consequence, if $\mathcal F$ is a regular foliation of rank $d\ge 3$ and $\wedge^2\mathcal F$ is strictly nef, then the fibers of the fibration in Theorem \ref{main} are projective spaces or smooth hyperquadrics. See Corollary \ref{coro}.

We now outline the proofs of the main results.

The proof of Theorem \ref{wedge3} starts with the observation that nefness of $\wedge^3T_X$ implies nefness of $-K_X$. Hence one may apply the decomposition theorem of Cao and H\"oring after passing to a finite \'etale cover. The variety then decomposes as a product $Y\times Z$, where $K_Y$ is trivial and the Albanese morphism of $Z$ is a locally trivial fibration with rationally connected fibers. The key point is to analyze the relative tangent sequence. The induced filtration on $\wedge^3T_X$ shows that, after restricting to fibers, either $T_F$ or $\wedge^2T_F$ is nef, depending on the dimension of the base. This reduces the proof to Watanabe's theorem for $\wedge^2T_X$, to the rationally connected case, and to a Calabi--Yau case handled by a Chern class computation together with the characterization of torus quotients in \cite{greb2016movable}.

For Theorem \ref{maintheorem}, we follow the strategy of Ross. Let \(\mathcal K\) be a minimal covering family of rational curves. There are two possible splitting types for \(f^*\mathcal E\) along curves in \(\mathcal K\). In the first case, a slope argument with respect to the curve class of \(\mathcal K\) shows that the rational quotient must be trivial, and hence \(X\simeq \PP^n\). In the second case, we argue by induction on the dimension. The rational quotient reduces the problem to a projective bundle over a curve. A vanishing result for such projective bundles then rules out the existence of the required subsheaf unless the Picard number is one; Ross's result applies in that case.

The proof of Theorem \ref{main} uses Ou's structure theorem for foliations with nef anti-canonical class. It gives a locally trivial fibration \(f:X\to Y\) and a transverse foliation \(\mathcal G\) on \(Y\) with \(K_{\mathcal G}\equiv 0\). The strict nefness of \(\wedge^r\mathcal F\) forces this transverse foliation to vanish, so \(\mathcal F=T_{X/Y}\). The fibers are then Fano by Watanabe's result on strictly nef exterior powers. Finally, the monodromy representation of the locally constant fibration is used, via Shafarevich maps, to show that \(Y\) is hyperbolic and that \(K_Y\) is ample.

\noindent\textbf{Acknowledgments}. The author wishes to express his gratitude to Jie Liu and Xiaokui Yang for suggesting the problems. He also would like to thank Jie Liu for many insightful discussions and for his careful reading of the manuscript. The author further appreciates the encouragement and helpful discussions provided by Junyu Cao, Paolo Cascini, Zhizhong Huang, Wenhao Ou, and Zhiwei Zheng.
\section{Preliminaries}\label{pre}

Throughout the paper, we work over the field $\mathbb C$.

\subsection{Positivity of divisors and vector bundles}

\bdefinition
Let $X$ be a projective variety. A Cartier divisor $D$ on $X$ is called
nef, respectively strictly nef, if $D\cdot C\ge 0$, respectively
$D\cdot C>0$, for every irreducible curve $C\subset X$.
\edefinition

\bdefinition
Let $E$ be a vector bundle on a projective variety $X$. We say that $E$
is nef, respectively strictly nef or ample, if the tautological line bundle
$\OO_{\PP(E)}(1)$ is nef, respectively strictly nef or ample, on
$\PP(E)=\mathbf{Proj}\,\operatorname{Sym}^{\bullet}E$.
\edefinition

The following standard properties will be used repeatedly.

\bproposition[\cite{lazarsfeld2017positivity}, Theorem 6.2.12 and
\cite{campana1991projective}, Proposition 1.2]\label{nef}
Let $0\to F\to E\to G\to 0$ be a short exact sequence of vector bundles
on a projective variety. Then the following hold.

\noindent\textup{(1)} Every quotient bundle of a nef vector bundle is nef.

\noindent\textup{(2)} If $E$ is nef, then $\wedge^kE$ is nef for every $k>0$.

\noindent\textup{(3)} If both $F$ and $G$ are nef, then $E$ is nef.

\noindent\textup{(4)} If $E$ is nef and $\det G$ is numerically trivial,
then $F$ is nef.
\eproposition

\subsection{Foliations and locally constant fibrations}\label{foliationpre}

\bdefinition
Let $X$ be a projective manifold. A foliation on $X$ is a coherent
subsheaf $\mathcal F\subset T_X$ which is closed under the Lie bracket and
saturated in $T_X$. It is called regular if it is a subbundle of $T_X$.
Its canonical class is defined by
$\OO_X(-K_{\mathcal F})\simeq\det\mathcal F$. The foliation is called
algebraically integrable if it is induced by a dominant rational map.
\edefinition

\bdefinition
Let $f:X\to Y$ be a surjective morphism with connected fibers between
normal varieties. We say that $f$ is a locally constant fibration with
fiber $F$ if it is analytically locally trivial and there exists a
representation $\rho:\pi_1(Y)\to\operatorname{Aut}(F)$ such that
\[
X\simeq (Y^{\mathrm{univ}}\times F)/\pi_1(Y),
\]
where $\gamma\in\pi_1(Y)$ acts by
$\gamma\cdot(y,z)=(\gamma\cdot y,\rho(\gamma)(z))$.
\edefinition

A vector bundle $E$ of rank $r$ on a normal variety $Y$ is called flat if
it is induced by a representation $\pi_1(Y)\to\operatorname{GL}(r,\mathbb C)$.

\section{Nefness of the third exterior power of tangent bundles}

In this section we prove Theorem \ref{wedge3}. The main input is the decomposition theorem of Cao and H\"oring for projective manifolds with nef anti-canonical bundle.

\btheorem[\cite{cao2017decomposition}, Theorem 1.2]\label{cao}
Let $X$ be a smooth projective variety with nef $-K_X$. Then there exists a finite \'etale cover $f:\tilde X\to X$ such that
\[
\tilde X\simeq Y\times Z,
\]
where $\omega_Y\simeq \OO_Y$, and the Albanese morphism $\alpha_Z:Z\to \operatorname{Alb}(Z)$ is a locally trivial fibration whose fibers are rationally connected.
\etheorem

We will also use the following results of Watanabe.

\btheorem[\cite{watanabe2021positivity}, Corollary 1.6]\label{wanneffano}
Let $X$ be a projective manifold such that $\wedge^2T_X$ is nef and $\dim X\ge 3$. Then either $T_X$ is nef or $X$ is a Fano manifold.
\etheorem

\btheorem[\cite{watanabe2022positivity}, Proposition 1.4]\label{RCfano}
Let $X$ be a projective manifold of dimension $n$. If $X$ is rationally connected and $\wedge^rT_X$ is nef for some $1\le r<n$, then $X$ is Fano.
\etheorem

We first record two elementary consequences of the nefness of $\wedge^3T_X$.

\blemma\label{KX}
Let $X$ be a smooth projective variety of dimension $n\ge 3$. If $\wedge^3T_X$ is nef, then $-K_X$ is nef.
\elemma

\bproof
Since
\[
\det(\wedge^3T_X)\simeq (-K_X)^{\otimes \binom{n-1}{2}},
\]
and the determinant of a nef vector bundle is nef, the divisor $-K_X$ is nef.
\eproof

\blemma\label{fibernef}
Let $X$ be a smooth projective variety with $\wedge^3T_X$ nef, and assume that $\dim X\ge 4$. Let $\varphi:X\to Y$ be a smooth morphism with irreducible fibers. Then the following statements hold.

\noindent\textup{(1)} If $\dim Y\ge 3$, then $\wedge^3T_Y$ is nef.

\noindent\textup{(2)} If $\dim Y\ge 2$, then every fiber $F$ of $\varphi$ has nef tangent bundle $T_F$.

\noindent\textup{(3)} If $\dim Y=1$, then every fiber $F$ of $\varphi$ has nef $\wedge^2T_F$.
\elemma

\bproof
Consider the relative tangent sequence
\[
0\longrightarrow T_{X/Y}\longrightarrow T_X\longrightarrow \varphi^*T_Y\longrightarrow 0.
\]
It induces a natural filtration of $\wedge^3T_X$ whose graded pieces are
\[
\wedge^3T_{X/Y},\qquad
\wedge^2T_{X/Y}\otimes \varphi^*T_Y,\qquad
T_{X/Y}\otimes \varphi^*\wedge^2T_Y,\qquad
\varphi^*\wedge^3T_Y.
\]

If $\dim Y\ge 3$, then $\varphi^*\wedge^3T_Y$ is a quotient of $\wedge^3T_X$, hence it is nef. Since $\varphi$ is surjective, this implies that $\wedge^3T_Y$ is nef.

Now let $F$ be a fiber of $\varphi$, and set $q=\dim Y$. Restricting the above filtration to $F$, its graded pieces become
\[
\wedge^3T_F,\qquad
(\wedge^2T_F)^{\oplus q},\qquad
T_F^{\oplus \binom{q}{2}},\qquad
\OO_F^{\oplus \binom{q}{3}}.
\]
Since $\wedge^3T_X|_F$ is nef and the last graded piece is trivial, Proposition \ref{nef} implies that the kernel of the quotient map to $\OO_F^{\oplus \binom{q}{3}}$ is nef. If $q\ge 2$, then $T_F^{\oplus \binom{q}{2}}$ is a quotient of this nef bundle, hence $T_F$ is nef. This proves \textup{(2)}.

If $q=1$, the filtration reduces to an exact sequence whose quotient is $\wedge^2T_F$. Therefore $\wedge^2T_F$ is nef, proving \textup{(3)}.
\eproof

The preceding lemma immediately gives the following product case.

\bcorollary\label{prod}
Let $X$ be a smooth projective variety such that $X\simeq Y\times Z$, where $\dim Y\ge 2$ and $\dim Z\ge 2$. If $\wedge^3T_X$ is nef, then $T_X$ is nef.
\ecorollary

\bproof
Apply Lemma \ref{fibernef} to the two projections $X\to Y$ and $X\to Z$. We obtain that both $T_Y$ and $T_Z$ are nef. Hence
\[
T_X\simeq p_Y^*T_Y\oplus p_Z^*T_Z
\]
is nef.
\eproof

We also need to treat the case where the Calabi--Yau factor is the whole finite \'etale cover.

\blemma\label{wxtrivial}
Let $X$ be a smooth projective variety of dimension $n\ge 4$. Assume that
$K_X\simeq\OO_X$ and that $\wedge^3T_X$ is nef. Then $X$ admits a finite
\'{e}tale cover by an abelian variety. In particular, $T_X$ is numerically
flat, and hence nef.
\elemma

\bproof
Since $K_X\simeq\OO_X$, one has $c_1(T_X)=0$ and therefore
$c_1(\wedge^3T_X)=0$. The nef bundle $\wedge^3T_X$ is thus numerically
flat. In particular,
\[
c_2(\wedge^3T_X)\cdot H^{n-2}=0
\]
for every ample divisor $H$ on $X$. By the splitting principle,
\[
c_2(\wedge^3T_X)
=
\frac{1}{2}\binom{n-1}{2}
\left(\binom{n-1}{2}-1\right)c_1^2(T_X)
+
\frac{(n-2)(n-3)}{2}c_2(T_X).
\]
Since $c_1(T_X)=0$ and $n\ge4$, it follows that
$c_2(X)\cdot H^{n-2}=0$. By
\cite[Theorem 7.1]{greb2016movable}, there exists an abelian variety $A$
and a finite quasi-\'{e}tale morphism $\gamma:A\to X$. Since $X$ is
smooth, purity of the branch locus implies that $\gamma$ is \'{e}tale.
Consequently,
$\gamma^*T_X\simeq T_A\simeq\OO_A^{\oplus n}$, so $T_X$ is numerically
flat.
\eproof

We can now prove the structure theorem.

\btheorem\label{wedge31}
Let $X$ be a smooth projective variety with $\wedge^3T_X$ nef and $\dim X\ge 4$. Then at least one of the following holds.

\noindent\textup{(1)} The manifold $X$ is Fano.

\noindent\textup{(2)} The tangent bundle $T_X$ is nef.

\noindent\textup{(3)} There exists a finite \'etale cover $\tilde X\to X$ such that $\tilde X$ admits a locally trivial Fano fibration over an elliptic curve.
\etheorem

\bproof
By Lemma \ref{KX}, the anti-canonical divisor $-K_X$ is nef. Applying Theorem \ref{cao}, we obtain a finite \'etale cover $f:\tilde X\to X$ such that
\[
\tilde X\simeq Y\times Z,
\]
where $\omega_Y\simeq \OO_Y$, and the Albanese morphism $\alpha_Z:Z\to \operatorname{Alb}(Z)$ is a locally trivial fibration whose fibers are rationally connected. Since $f$ is finite \'etale, the bundle $\wedge^3T_{\tilde X}$ is nef.

It is enough to prove the desired alternatives for $\tilde X$. Indeed, nefness of the tangent bundle descends under finite \'etale covers, since $T_{\tilde X}\simeq f^*T_X$. Similarly, if $\tilde X$ is Fano, then $X$ is Fano.

We distinguish several cases according to the dimensions of $Y$ and $Z$.

\smallskip
\noindent\emph{Case 1.} Assume that $\dim Y\ge 2$ and $\dim Z\ge 2$.
In this case, Corollary \ref{prod} shows that $T_{\tilde X}$ is nef. Hence $T_X$ is nef.

\smallskip
\noindent\emph{Case 2.} Assume that $\dim Y=1$.
Since $\omega_Y\simeq \OO_Y$, the curve $Y$ is elliptic. The assumption $\dim X\ge 4$ implies $\dim Z\ge 3$. Applying Lemma \ref{fibernef} to the projection $\tilde X=Y\times Z\to Y$, we see that $\wedge^2T_Z$ is nef. By Theorem \ref{wanneffano}, either $T_Z$ is nef or $Z$ is Fano.

If $T_Z$ is nef, then $T_{\tilde X}$ is nef, hence $T_X$ is nef. If $Z$ is Fano, then the projection $\tilde X=Y\times Z\to Y$ is a locally trivial Fano fibration over the elliptic curve $Y$.

\smallskip
\noindent\emph{Case 3.} Assume that $\dim Z=1$.
The Albanese morphism $\alpha_Z:Z\to \operatorname{Alb}(Z)$ has rationally connected fibers.

If $\dim \operatorname{Alb}(Z)=0$, then $Z\simeq \PP^1$. Applying Lemma \ref{fibernef} to the projection $\tilde X=\PP^1\times Y\to \PP^1$, we obtain that $\wedge^2T_Y$ is nef. By Theorem \ref{wanneffano}, either $T_Y$ is nef or $Y$ is Fano. In the first case, $T_{\tilde X}$ is nef. In the second case, $\tilde X$ is Fano. Thus the conclusion follows.

If $\dim \operatorname{Alb}(Z)=1$, then $Z$ is an elliptic curve. Again, applying Lemma \ref{fibernef} to the projection $\tilde X=Y\times Z\to Z$, we obtain that $\wedge^2T_Y$ is nef. By Theorem \ref{wanneffano}, either $T_Y$ is nef or $Y$ is Fano. In the first case, $T_{\tilde X}$ is nef. In the second case, $\tilde X\to Z$ is a locally trivial Fano fibration over an elliptic curve.

\smallskip
\noindent\emph{Case 4.} Assume that $\dim Z=0$.
Then $\tilde X\simeq Y$ and $K_{\tilde X}\simeq \OO_{\tilde X}$. By Lemma \ref{wxtrivial}, the tangent bundle $T_{\tilde X}$ is nef. Hence $T_X$ is nef.

\smallskip
\noindent\emph{Case 5.} Assume that $\dim Y=0$.
Then $\tilde X\simeq Z$, and the Albanese morphism
\[
\alpha:\tilde X\to A:=\operatorname{Alb}(\tilde X)
\]
is a locally trivial fibration with rationally connected fibers.

If $\dim A=0$, then $\tilde X$ is rationally connected. Since $\wedge^3T_{\tilde X}$ is nef and $3<\dim \tilde X$, Theorem \ref{RCfano} implies that $\tilde X$ is Fano. Hence $X$ is Fano.

Assume next that $\dim A=q\ge 2$. Since $T_A$ is trivial, the relative tangent sequence
\[
0\longrightarrow T_{\tilde X/A}\longrightarrow T_{\tilde X}\longrightarrow \alpha^*T_A\longrightarrow 0
\]
induces a filtration of $\wedge^3T_{\tilde X}$ whose graded pieces are
\[
\wedge^3T_{\tilde X/A},\qquad
(\wedge^2T_{\tilde X/A})^{\oplus q},\qquad
T_{\tilde X/A}^{\oplus \binom{q}{2}},\qquad
\OO_{\tilde X}^{\oplus \binom{q}{3}}.
\]
As before, the kernel of the quotient map to the trivial last graded piece is nef by Proposition \ref{nef}. Hence $T_{\tilde X/A}^{\oplus \binom{q}{2}}$ is nef, and so $T_{\tilde X/A}$ is nef. Since $\alpha^*T_A$ is trivial, Proposition \ref{nef} applied to the relative tangent sequence shows that $T_{\tilde X}$ is nef. Therefore $T_X$ is nef.

Finally, assume that $\dim A=1$. Then $A$ is an elliptic curve. Applying Lemma \ref{fibernef} to $\alpha:\tilde X\to A$, we obtain that $\wedge^2T_F$ is nef for every fiber $F$ of $\alpha$. The fibers of $\alpha$ are rationally connected by Theorem \ref{cao}, and $\dim F\ge 3$. Therefore Theorem \ref{RCfano} implies that every fiber $F$ is Fano. Hence $\alpha:\tilde X\to A$ is a locally trivial Fano fibration over an elliptic curve. This completes the proof.\eproof

\section{Ample subsheaves of the second exterior power of the tangent bundle}

In this section, we prove Theorem \ref{maintheorem}. We begin with two examples showing that the hypothesis
\(\mathcal F=\wedge^p\mathcal E\) in Conjecture \ref{conjecture} cannot be omitted. We then recall the results on minimal rational curves and rational quotients that will be used throughout the proof. By Lemma \ref{twocase}, there are two possible splitting types for \(\mathcal E\) along a minimal covering family. The non-uniform splitting type is treated by a slope argument, whereas the uniform splitting type is handled by induction on the dimension. The latter argument eventually reduces the problem to a projective bundle over a curve.

\bexample
Let \(X=\PP^m\times \PP^n\), and let \(p_1\) and \(p_2\) be the two projections. Since
\[
T_X\simeq p_1^*T_{\PP^m}\oplus p_2^*T_{\PP^n},
\]
the bundle \(p_1^*T_{\PP^m}\otimes p_2^*T_{\PP^n}\) is a direct summand of \(\wedge^2T_X\). To see that it is ample, recall from the Euler sequences that
\(T_{\PP^m}\) and \(T_{\PP^n}\) are quotients of
\(\OO_{\PP^m}(1)^{\oplus(m+1)}\) and
\(\OO_{\PP^n}(1)^{\oplus(n+1)}\), respectively. Hence
\(p_1^*T_{\PP^m}\otimes p_2^*T_{\PP^n}\) is a quotient of
\[
\OO_X(1,1)^{\oplus(m+1)(n+1)},
\]
which is ample. In general, this subsheaf is not of the form
\(\wedge^2\mathcal E\).
\eexample

\bexample
Let \(Y\) be a projective manifold, and let \(L_1,\ldots,L_{p+1}\) be ample line bundles on \(Y\) such that
\[
H^0\left(Y,L_1^{\otimes(p-1)}\otimes
(L_2\otimes\cdots\otimes L_{p+1})^{-1}\right)\neq 0.
\]
Set
\[
\mathcal V:=\bigoplus_{i=1}^{p+1}L_i
\qquad\text{and}\qquad
X:=\PP_Y(\mathcal V).
\]
Then \(\OO_{\PP(\mathcal V)}(1)\) is ample and there is an inclusion
\[
\OO_{\PP(\mathcal V)}(1)\subset \wedge^pT_X.
\]
Indeed,
\[
\begin{aligned}
H^0\left(X,\wedge^pT_{X/Y}\otimes\OO_{\PP(\mathcal V)}(-1)\right)
&\simeq H^0\left(X,\pi^*\det(\mathcal V)^*\otimes
\OO_{\PP(\mathcal V)}(p)\right)\\
&\simeq H^0\left(Y,\det(\mathcal V)^*\otimes S^p\mathcal V\right).
\end{aligned}
\]
The line bundle
\[
L_1^{\otimes(p-1)}\otimes
(L_2\otimes\cdots\otimes L_{p+1})^{-1}
\]
is a direct summand of \(\det(\mathcal V)^*\otimes S^p\mathcal V\), so the assumed section gives the desired inclusion.
\eexample

\subsection{Rational quotients and splitting types}\label{sec:rq-splitting}

We first recall the structure theorem for rational quotients that underlies both parts of the proof.

\btheorem[\cite{araujo2007cohomological}, Theorem 2.7]\label{summary}
Let \(X\) be a smooth complex projective variety, let \(\mathcal K\) be a minimal covering family of rational curves on \(X\), and let
\[
\pi^\circ:X^\circ\to Y^\circ
\]
be the \(\mathcal K\)-rationally connected quotient. Suppose that \(T_X\) contains a subsheaf \(\mathcal D\) such that \(f^*\mathcal D\) is ample for a general member \([f]\in\mathcal K\). Then, after shrinking \(X^\circ\) and \(Y^\circ\) if necessary, the morphism \(\pi^\circ\) is a projective space bundle, and the inclusion
\(\mathcal D|_{X^\circ}\hookrightarrow T_{X^\circ}\) factors through the natural inclusion
\(T_{X^\circ/Y^\circ}\hookrightarrow T_{X^\circ}\). Moreover, if \(\mathcal K\) is unsplit, then the rational quotient extends to an open subset of \(X\) whose complement has codimension at least two.
\etheorem

We also recall two results of Ross.

\blemma[\cite{ross2010characterizations}, Lemma 3.1.1]\label{rossuniruled}
Let \(X\) be a smooth complex projective variety, let \(\mathcal E\) be an ample vector bundle of rank \(r\) on \(X\), and suppose that
\(\wedge^p\mathcal E\subseteq\wedge^pT_X\) for some positive integer \(p\le r\). Then \(X\) is uniruled. In particular, \(X\) admits a minimal covering family of rational curves.
\elemma

\blemma[\cite{ross2010characterizations}, Lemma 3.1.2 and Corollary 3.1.3]\label{twocase}
Let \(X\) be a smooth complex projective variety, let \(\mathcal E\) be an ample vector bundle of rank \(r\) on \(X\), and suppose that
\(\wedge^p\mathcal E\subseteq\wedge^pT_X\) for some positive integer \(p\le r\). Let \(\mathcal K\) be a minimal covering family of rational curves on \(X\). Then one of the following alternatives holds for every member \(f:\PP^1\to X\) of \(\mathcal K\):
\[
f^*\mathcal E\simeq
\OO_{\PP^1}(2)\oplus\OO_{\PP^1}(1)^{\oplus r-1},
\]
or
\[
f^*\mathcal E\simeq\OO_{\PP^1}(1)^{\oplus r}.
\]
Moreover, the family \(\mathcal K\) is unsplit.
\elemma

The following elementary lemma will be used to show that a morphism to \(\wedge^2T_X\) is vertical with respect to a rational quotient.

\blemma\label{shortexactsequence}
Let \(\pi:X\to Y\) be a smooth morphism between smooth varieties, and let \(\mathcal V\) be a vector bundle on \(X\). Assume that \(\mathcal V|_F\) is ample for a general fiber \(F\) of \(\pi\). Then there is an exact sequence
\[
\begin{aligned}
0\longrightarrow
&H^0\left(X,\wedge^2T_{X/Y}\otimes\mathcal V^*\right)
\longrightarrow H^0\left(X,\wedge^2T_X\otimes\mathcal V^*\right)\\
&\longrightarrow
H^0\left(X,T_{X/Y}\otimes\pi^*T_Y\otimes\mathcal V^*\right).
\end{aligned}
\]
\elemma

\bproof
The relative tangent sequence
\[
0\longrightarrow T_{X/Y}\longrightarrow T_X\longrightarrow\pi^*T_Y\longrightarrow 0
\]
induces a filtration
\[
0=\mathcal G_3\subset\mathcal G_2\subset\mathcal G_1\subset
\mathcal G_0=\wedge^2T_X\otimes\mathcal V^*
\]
whose graded pieces are
\[
\mathcal G_i/\mathcal G_{i+1}
\simeq
\wedge^iT_{X/Y}\otimes\pi^*\wedge^{2-i}T_Y\otimes\mathcal V^*
\]
for \(0\le i\le 2\). Since \(\mathcal V|_F\) is ample for a general fiber \(F\), one has
\[
H^0\left(X,\pi^*\wedge^2T_Y\otimes\mathcal V^*\right)=0.
\]
Hence
\[
H^0\left(X,\wedge^2T_X\otimes\mathcal V^*\right)
=H^0(X,\mathcal G_1),
\]
and the assertion follows from the exact sequence
\[
0\longrightarrow\wedge^2T_{X/Y}\otimes\mathcal V^*
\longrightarrow\mathcal G_1
\longrightarrow T_{X/Y}\otimes\pi^*T_Y\otimes\mathcal V^*
\longrightarrow 0.
\]
\eproof

By Lemma \ref{twocase}, the proof of Theorem \ref{maintheorem} naturally separates according to the two splitting types above. We treat the non-uniform splitting type in Subsection \ref{sec:nonuniform}. The uniform splitting type will be handled by induction after establishing, in Subsection \ref{sec:projective-curves}, a vanishing result for projective bundles over curves.

\subsection{The non-uniform splitting case}\label{sec:nonuniform}

Let \(X\) be a projective manifold of dimension \(n\), and let \(\mathcal E\) be an ample vector bundle of rank \(r\ge2\) such that
\[
\wedge^2\mathcal E\subset\wedge^2T_X.
\]
Let \(\mathcal K\) be an unsplit minimal covering family of rational curves on \(X\). In this subsection, we assume that, for a general member \(f:\PP^1\to X\) of \(\mathcal K\),
\[
f^*T_X\simeq
\OO_{\PP^1}(2)\oplus
\OO_{\PP^1}(1)^{\oplus d}\oplus
\OO_{\PP^1}^{\oplus n-d-1}
\]
and
\[
f^*\mathcal E\simeq
\OO_{\PP^1}(2)\oplus
\OO_{\PP^1}(1)^{\oplus r-1}.
\]
We prove that this splitting type already forces \(X\simeq\PP^n\).

\blemma\label{2.1}
Under the assumptions above, one has \(r-1\le d\).
\elemma

\bproof
Restrict the inclusion \(\wedge^2\mathcal E\subset\wedge^2T_X\) to a general curve in \(\mathcal K\). The summand \(\OO_{\PP^1}(3)\) occurs in \(\wedge^2f^*\mathcal E\) with multiplicity \(r-1\), and in \(\wedge^2f^*T_X\) with multiplicity \(d\). Therefore \(r-1\le d\).
\eproof

We recall the slope notation only where it is needed. Let $\alpha$ be a
movable curve class and let $\mathcal G$ be a torsion-free sheaf. Set
\[
\mu_\alpha(\mathcal G)
=
\frac{c_1(\mathcal G)\cdot\alpha}{\operatorname{rank}\mathcal G},
\qquad
\mu_\alpha^{\max}(\mathcal G)
=
\sup_{0\neq\mathcal A\subset\mathcal G}\mu_\alpha(\mathcal A).
\]
The maximal slope is attained by a unique maximal destabilizing subsheaf;
see \cite[Corollary 2.25]{greb2016movable}. We say that $\mathcal G$ is
$\alpha$-semistable if
$\mu_\alpha^{\max}(\mathcal G)=\mu_\alpha(\mathcal G)$.

\blemma\label{2.2}
Under the assumptions at the beginning of this subsection, let
$\alpha=[C]\in N_1(X)$ be the numerical class of a curve parametrized by
$\mathcal K$. Then
\[
\mu_\alpha^{\max}(T_X)\ge\mu_\alpha(\mathcal E)>1.
\]
\elemma

\bproof
By \cite[Theorem 4.2]{greb2016movable}, one has the equality
\[
\mu_\alpha^{\max}(T_X\otimes T_X)
=
2\mu_\alpha^{\max}(T_X).
\]
Since $\wedge^2T_X$ is a direct summand of $T_X\otimes T_X$ and
$\wedge^2\mathcal E\subset\wedge^2T_X$, it follows that
\[
2\mu_\alpha^{\max}(T_X)
\ge
\mu_\alpha^{\max}(\wedge^2T_X)
\ge
\mu_\alpha(\wedge^2\mathcal E)
=
2\mu_\alpha(\mathcal E).
\]
Finally, the splitting
$ f^*\mathcal E\simeq\OO_{\PP^1}(2)\oplus
\OO_{\PP^1}(1)^{\oplus r-1}$ gives
$\mu_\alpha(\mathcal E)=(r+1)/r>1$.
\eproof

\btheorem\label{firstcase}
Let \(X\) be a projective manifold of dimension \(n\), and let \(\mathcal E\) be an ample vector bundle of rank \(r\ge2\) such that
\(\wedge^2\mathcal E\subset\wedge^2T_X\). Let \(\mathcal K\) be an unsplit minimal covering family of rational curves on \(X\). Assume that, for a general member \(f:\PP^1\to X\) of \(\mathcal K\),
\[
f^*T_X\simeq
\OO_{\PP^1}(2)\oplus
\OO_{\PP^1}(1)^{\oplus d}\oplus
\OO_{\PP^1}^{\oplus n-d-1}
\]
and
\[
f^*\mathcal E\simeq
\OO_{\PP^1}(2)\oplus
\OO_{\PP^1}(1)^{\oplus r-1}.
\]
Then \(X\simeq\PP^n\).
\etheorem

\bproof
Let \(\alpha\) be the numerical class of a curve in \(\mathcal K\). Suppose first that \(T_X\) is \(\alpha\)-semistable. By Lemma \ref{2.2},
\[
\mu_\alpha(T_X)=\mu_\alpha^{\max}(T_X)>1.
\]
On the other hand,
\[
\mu_\alpha(T_X)=\frac{d+2}{n}.
\]
Since \(d\le n-1\), it follows that \(d=n-1\). Thus
\[
f^*T_X\simeq
\OO_{\PP^1}(2)\oplus\OO_{\PP^1}(1)^{\oplus n-1}.
\]
Applying Theorem \ref{summary} with \(\mathcal D=T_X\), we conclude that the rational quotient has zero-dimensional base and that \(X\simeq\PP^n\).

Assume now that \(T_X\) is not \(\alpha\)-semistable, and let \(\mathcal D\subset T_X\) be its maximal destabilizing subsheaf. By Lemma \ref{2.2},
\(\mu_\alpha(\mathcal D)>1\). For a general member \(f:\PP^1\to X\) of \(\mathcal K\), this implies
\[
f^*\mathcal D\simeq
\OO_{\PP^1}(2)\oplus
\OO_{\PP^1}(1)^{\oplus r'-1},
\]
where \(r'=\operatorname{rank}\mathcal D\). By Theorem \ref{summary}, the rational quotient
\[
\pi:X^\circ\to Y^\circ
\]
is a \(\PP^k\)-bundle on an open subset \(X^\circ\subset X\) whose complement has codimension at least two. Here \(k=d+1\), and the inclusion
\(\mathcal D|_{X^\circ}\subset T_{X^\circ}\) factors through
\(T_{X^\circ/Y^\circ}\).

We show that \(Y^\circ\) is a point. Suppose otherwise. We first prove that
\[
\wedge^2\mathcal E|_{X^\circ}
\subset
\wedge^2T_{X^\circ/Y^\circ}.
\]
Apply Lemma \ref{shortexactsequence} with
\(\mathcal V=\wedge^2\mathcal E|_{X^\circ}\). It is enough to prove
\[
H^0\left(X^\circ,
T_{X^\circ/Y^\circ}\otimes\pi^*T_{Y^\circ}
\otimes\wedge^2\mathcal E^*|_{X^\circ}\right)=0.
\]
If this group were nonzero, then, after restricting to a general fiber
\(F\simeq\PP^k\), we would obtain a nonzero morphism
\[
\wedge^2\mathcal E|_F\longrightarrow T_F.
\]
By \cite{liu2019characterization}, Theorem 1.1, its image is either
\(T_{\PP^k}\) or a direct sum of copies of \(\OO_{\PP^k}(1)\). Restricting to a general line \(\ell\subset F\), the source has splitting
\[
\wedge^2\mathcal E|_\ell\simeq
\OO_{\PP^1}(3)^{\oplus r-1}
\oplus
\OO_{\PP^1}(2)^{\oplus\binom{r-1}{2}}.
\]
It cannot have a nonzero quotient isomorphic to a direct sum of copies of
\(\OO_{\PP^1}(1)\). If the image is \(T_{\PP^k}\), then the restriction
\[
\wedge^2\mathcal E|_\ell
\longrightarrow
\OO_{\PP^1}(2)\oplus\OO_{\PP^1}(1)^{\oplus k-1}
\]
forces \(k=1\). But then Lemma \ref{2.1} gives
\(r-1\le d=k-1=0\), contradicting \(r\ge2\). The required vanishing follows, and hence the inclusion is vertical.

By Lemma \ref{2.1}, one has \(r\le k\). Suppose first that \(r<k\). Along a line parametrized by \(\mathcal K\),
\[
\mu_\alpha(\wedge^2\mathcal E)
=2\frac{r+1}{r}
>
2\frac{k+1}{k}
=\mu_\alpha(\wedge^2T_{X^\circ/Y^\circ}).
\]
Thus \(\wedge^2T_{X^\circ/Y^\circ}\) is not \(\alpha\)-semistable. Since exterior powers of semistable sheaves are semistable in characteristic zero,
\(T_{X^\circ/Y^\circ}\) is not \(\alpha\)-semistable either. Let
\(\mathcal D'\subset T_{X^\circ/Y^\circ}\) be a maximal destabilizing subsheaf. For a general line \(\ell\) in a general fiber of \(\pi\),
\[
\mathcal D'|_\ell\simeq
\OO_{\PP^1}(2)\oplus
\OO_{\PP^1}(1)^{\oplus s-1}.
\]
By \cite{andreatta2001manifolds}, Lemma 2.2, applied to
\(\mathcal D'\subset T_{X^\circ}\), one has
\[
(f^*T_{X^\circ})^+\subset f^*\mathcal D'
\]
for a general parametrization \(f:\PP^1\to\ell\). Since
\(\mathcal D'\subset T_{X^\circ/Y^\circ}\), the reverse inclusion holds, and therefore
\[
f^*\mathcal D'=(f^*T_{X^\circ})^+=f^*T_{X^\circ/Y^\circ}.
\]
It follows that \(\mathcal D'\) has the same rank as
\(T_{X^\circ/Y^\circ}\). Since \(\mathcal D'\) is saturated, one obtains
\(\mathcal D'=T_{X^\circ/Y^\circ}\), contradicting the fact that
\(\mathcal D'\) destabilizes \(T_{X^\circ/Y^\circ}\).

It remains to consider \(r=k\). Choose a complete curve
\(B\subset Y^\circ\) through a general point, form the base change, and suppress the pullback notation. We are reduced to a smooth \(\PP^k\)-bundle
\[
\pi:X\to B
\]
together with an inclusion
\[
\wedge^2\mathcal E\subset\wedge^2T_{X/B}.
\]
Since \(r=k\), taking determinants gives an inclusion
\[
(\det\mathcal E)^{\otimes(k-1)}
\hookrightarrow
(-K_{X/B})^{\otimes(k-1)}.
\]
On every fiber \(F\simeq\PP^k\), both line bundles restrict to
\(\OO_{\PP^k}((k-1)(k+1))\). Hence there exists an effective divisor \(D\), supported over finitely many points of \(B\), such that
\[
-(k-1)K_{X/B}
\sim
(k-1)\det\mathcal E+D.
\]
The pair \((X,\frac{1}{k-1}D)\) is log canonical over the generic point of \(B\), and
\[
-K_{X/B}-\frac{1}{k-1}D\sim\det\mathcal E
\]
is ample. This contradicts \cite{araujo2007cohomological}, Lemma 3.1.

Therefore \(\dim Y^\circ=0\), and Theorem \ref{summary} yields
\(X\simeq\PP^n\).
\eproof

\subsection{Projective bundles over curves}\label{sec:projective-curves}

The induction argument for the uniform splitting case reduces the problem to a projective bundle over a curve. We establish here the two facts needed in that reduction.

\blemma\label{curveamplecriterion}
Let \(C\) be a smooth projective curve, let \(G\) and \(M\) be vector bundles on \(C\), and let
\[
\pi:X=\PP_C(G)\to C
\]
be the natural projection. If
\[
\pi^*M\otimes\OO_{\PP(G)}(1)
\]
is ample on \(X\), then \(G\otimes M\) is ample on \(C\).
\elemma

\bproof
Let \(\mu_{\min}\) denote the minimal slope of a vector bundle on \(C\). Let
\[
G\twoheadrightarrow Q_G
\qquad\text{and}\qquad
M\twoheadrightarrow Q_M
\]
be the last quotients in the Harder--Narasimhan filtrations of \(G\) and \(M\), respectively. Then \(Q_G\) and \(Q_M\) are semistable and
\[
\mu(Q_G)=\mu_{\min}(G),
\qquad
\mu(Q_M)=\mu_{\min}(M).
\]

The quotient \(G\twoheadrightarrow Q_G\) induces a closed immersion
\[
\PP_C(Q_G)\hookrightarrow\PP_C(G).
\]
Writing \(\pi_Q:\PP_C(Q_G)\to C\) for the projection, the vector bundle
\[
\pi_Q^*M\otimes\OO_{\PP(Q_G)}(1)
\]
is ample. The quotient \(M\twoheadrightarrow Q_M\) therefore gives an ample quotient bundle
\[
\mathcal A_Q:=
\pi_Q^*Q_M\otimes\OO_{\PP(Q_G)}(1).
\]
Set
\[
a:=\operatorname{rank}Q_G,
\qquad
b:=\operatorname{rank}Q_M,
\qquad
\xi:=c_1\bigl(\OO_{\PP(Q_G)}(1)\bigr).
\]
Since \(\det\mathcal A_Q\) is ample on the \(a\)-dimensional variety
\(\PP_C(Q_G)\), one has
\[
c_1(\det\mathcal A_Q)^a>0.
\]
Now
\[
\det\mathcal A_Q
\simeq
\pi_Q^*(\det Q_M)\otimes\OO_{\PP(Q_G)}(b),
\]
and the projective bundle formula gives
\[
\begin{aligned}
c_1(\det\mathcal A_Q)^a
&=
\left(b\xi+\pi_Q^*c_1(Q_M)\right)^a\\
&=b^a\deg Q_G+ab^{a-1}\deg Q_M\\
&=ab^a\bigl(\mu(Q_G)+\mu(Q_M)\bigr).
\end{aligned}
\]
Consequently,
\[
\mu_{\min}(G)+\mu_{\min}(M)>0.
\]
The tensor products of the semistable graded pieces in the Harder--Narasimhan filtrations of \(G\) and \(M\) are semistable, and hence
\[
\mu_{\min}(G\otimes M)
=
\mu_{\min}(G)+\mu_{\min}(M)>0.
\]
A vector bundle on a smooth projective curve is ample if and only if its minimal slope is positive. Therefore \(G\otimes M\) is ample.
\eproof

\blemma\label{vanish}
Let \(\pi:X=\PP_C(G)\to C\) be a projective bundle over a smooth projective curve \(C\), and let \(M\) be a vector bundle on \(C\). Assume that
\[
\pi^*M\otimes\OO_{\PP(G)}(1)
\]
is ample. Then
\[
H^0\left(C,\wedge^2G^*\otimes\wedge^2M^*\right)=0.
\]
Consequently,
\[
H^0\left(X,
\wedge^2T_{X/C}\otimes
\pi^*\wedge^2M^*\otimes
\OO_{\PP(G)}(-2)
\right)=0.
\]
\elemma

\bproof
If either \(\operatorname{rank}G<2\) or \(\operatorname{rank}M<2\), the first assertion is immediate. We may therefore assume that both ranks are at least two.

By Lemma \ref{curveamplecriterion}, the vector bundle \(G\otimes M\) is ample on \(C\). Hence \((G\otimes M)^{\otimes2}\) is ample. Since
\(\wedge^2G\otimes\wedge^2M\) is a quotient of
\((G\otimes M)^{\otimes2}\), it is ample as well. Its dual has no nonzero global sections, and therefore
\[
H^0\left(C,\wedge^2G^*\otimes\wedge^2M^*\right)=0.
\]

We now prove the second assertion. If the relative dimension of \(\pi\) is one, then \(T_{X/C}\) has rank one and \(\wedge^2T_{X/C}=0\). Thus we may assume that the relative dimension is at least two.

The relative Euler sequence is
\[
0\longrightarrow\OO_X
\longrightarrow\pi^*G^*\otimes\OO_X(1)
\longrightarrow T_{X/C}
\longrightarrow0.
\]
Taking second exterior powers gives
\[
0\longrightarrow T_{X/C}
\longrightarrow\wedge^2\bigl(\pi^*G^*\otimes\OO_X(1)\bigr)
\longrightarrow\wedge^2T_{X/C}
\longrightarrow0.
\]
After tensoring by
\(\pi^*\wedge^2M^*\otimes\OO_X(-2)\), we obtain
\[
\begin{aligned}
0\longrightarrow&
T_{X/C}\otimes\pi^*\wedge^2M^*\otimes\OO_X(-2)\\
\longrightarrow&
\pi^*(\wedge^2G^*\otimes\wedge^2M^*)\\
\longrightarrow&
\wedge^2T_{X/C}\otimes\pi^*\wedge^2M^*\otimes\OO_X(-2)
\longrightarrow0.
\end{aligned}
\]
The first part gives
\[
H^0\left(X,\pi^*(\wedge^2G^*\otimes\wedge^2M^*)\right)=0.
\]
It is therefore enough to prove
\[
H^1\left(X,
T_{X/C}\otimes\pi^*\wedge^2M^*\otimes\OO_X(-2)
\right)=0.
\]

Tensoring the relative Euler sequence by
\(\pi^*\wedge^2M^*\otimes\OO_X(-2)\) gives
\[
\begin{aligned}
0\longrightarrow&
\pi^*\wedge^2M^*\otimes\OO_X(-2)\\
\longrightarrow&
\pi^*(G^*\otimes\wedge^2M^*)\otimes\OO_X(-1)\\
\longrightarrow&
T_{X/C}\otimes\pi^*\wedge^2M^*\otimes\OO_X(-2)
\longrightarrow0.
\end{aligned}
\]
Since \(R^i\pi_*\OO_X(-1)=0\) for all \(i\), one has
\[
H^1\left(X,
\pi^*(G^*\otimes\wedge^2M^*)\otimes\OO_X(-1)
\right)=0.
\]
Moreover, because the relative dimension of \(\pi\) is at least two,
\[
R^i\pi_*\OO_X(-2)=0
\]
for all \(i\). Hence
\[
H^2\left(X,
\pi^*\wedge^2M^*\otimes\OO_X(-2)
\right)=0.
\]
The long exact cohomology sequence now gives the desired vanishing.
\eproof

\subsection{The uniform splitting case and the proof of Theorem B}\label{sec:uniform}

We now treat the second alternative in Lemma \ref{twocase}. For the induction on the dimension, it is convenient to prove the following slightly stronger statement.

\btheorem[A stronger form of Theorem \ref{maintheorem}]\label{rigid}
Let \(X\) be a projective manifold of dimension \(n\ge2\), and let \(\mathcal F\) be an ample locally free subsheaf of \(\wedge^2T_X\). Assume that
\(\mathcal F=\wedge^2\mathcal E\) for some ample vector bundle \(\mathcal E\). If
\[
\operatorname{rank}\mathcal F
\ge
\min\left\{\binom{n}{2},n\right\},
\]
then either \(X\simeq\PP^n\), or \(n=2\) and \(X\simeq Q_2\).
\etheorem

\bremark
When \(n\ge3\), one has \(\binom{n}{2}\ge n\), so Theorem \ref{rigid} implies Theorem \ref{maintheorem}. When \(n=2\), the rank condition is automatic, and the assertion is known by \cite{ross2010characterizations}.
\eremark

\bproof
We argue by induction on \(n\). The case \(n=2\) follows from \cite{ross2010characterizations}. Assume that the theorem holds in dimensions at most \(n-1\), where \(n\ge3\).

Let \(X\) be a projective manifold of dimension \(n\), and let
\[
\mathcal F=\wedge^2\mathcal E\subset\wedge^2T_X
\]
satisfy the assumptions of the theorem. By Lemma \ref{rossuniruled}, the manifold \(X\) is uniruled. Let \(\mathcal K\) be a minimal covering family of rational curves on \(X\). The non-uniform splitting type in Lemma \ref{twocase} is treated by Theorem \ref{firstcase}. We may therefore assume that
\[
f^*\mathcal E\simeq\OO_{\PP^1}(1)^{\oplus r}
\]
for every member \(f:\PP^1\to X\) of \(\mathcal K\), where
\(r=\operatorname{rank}\mathcal E\).

If \(r=2\), then \(\mathcal F=\det\mathcal E\) is an ample line bundle contained in \(\wedge^2T_X\), and the conclusion follows from \cite{druel2013characterizations}. Thus we may assume that \(r\ge3\).

If \(\rho(X)=1\), the theorem follows from \cite{ross2010characterizations}. Suppose, for contradiction, that \(\rho(X)>1\). Let
\[
\pi:X^\circ\to Y^\circ
\]
be the rational quotient of \(X\) with respect to \(\mathcal K\). Since \(\mathcal K\) is unsplit and \(\rho(X)>1\), the base \(Y^\circ\) has positive dimension. After shrinking \(X^\circ\) and \(Y^\circ\), we may assume that \(\pi\) is smooth.

We first prove that
\[
\wedge^2\mathcal E|_{X^\circ}
\subset
\wedge^2T_{X^\circ/Y^\circ}.
\]
Applying Lemma \ref{shortexactsequence} with
\(\mathcal V=\wedge^2\mathcal E|_{X^\circ}\), it is enough to show that
\[
H^0\left(X^\circ,
T_{X^\circ/Y^\circ}\otimes\pi^*T_{Y^\circ}
\otimes\wedge^2\mathcal E^*|_{X^\circ}
\right)=0.
\]
Suppose otherwise. Restricting a nonzero section to a general fiber \(F_y\) of \(\pi\), we obtain a nonzero morphism
\[
\wedge^2\mathcal E|_{F_y}\longrightarrow T_{F_y}.
\]
By \cite{liu2019characterization}, Theorem 1.1, the fiber \(F_y\) is isomorphic to \(\PP^k\), and the image is either \(T_{\PP^k}\) or a direct sum of copies of \(\OO_{\PP^k}(1)\). For a line \(\ell\subset F_y\), one has
\[
\wedge^2\mathcal E|_\ell
\simeq
\OO_{\PP^1}(2)^{\oplus\binom{r}{2}}.
\]
Thus the image cannot contain a summand \(\OO_{\PP^1}(1)\). If the image is \(T_{\PP^k}\), the splitting
\[
T_{\PP^k}|_\ell
\simeq
\OO_{\PP^1}(2)\oplus\OO_{\PP^1}(1)^{\oplus k-1}
\]
forces \(k=1\).

It remains to exclude this last possibility. If the general fiber is \(\PP^1\), then the relative tangent sequence induces an exact sequence
\[
0\longrightarrow
T_{X^\circ/Y^\circ}\otimes\pi^*T_{Y^\circ}
\longrightarrow
\wedge^2T_{X^\circ}
\longrightarrow
\pi^*\wedge^2T_{Y^\circ}
\longrightarrow0.
\]
The restriction of \(\wedge^2\mathcal E\) to a fiber is a direct sum of copies of \(\OO_{\PP^1}(2)\), whereas the restriction of
\(\pi^*\wedge^2T_{Y^\circ}\) is trivial. Hence the image of
\(\wedge^2\mathcal E\) in \(\pi^*\wedge^2T_{Y^\circ}\) is zero, and
\[
\wedge^2\mathcal E|_{X^\circ}
\subset
T_{X^\circ/Y^\circ}\otimes\pi^*T_{Y^\circ}.
\]
The right-hand side has rank \(n-1\), whereas
\(\operatorname{rank}\wedge^2\mathcal E=\operatorname{rank}\mathcal F\ge n\). This is impossible. The desired factorization follows.

Let \(F_y\) be a general fiber of \(\pi\), and set \(k=\dim F_y<n\). Restricting the factorization to \(F_y\), we obtain
\[
\wedge^2\mathcal E|_{F_y}\subset\wedge^2T_{F_y}.
\]
The bundle \(\wedge^2\mathcal E|_{F_y}\) is ample and satisfies the rank condition required by the induction hypothesis. Therefore \(F_y\) is either \(\PP^k\), or \(k=2\) and \(F_y\simeq Q_2\). The latter case is impossible, since \(\wedge^2T_{F_y}\) has rank one, whereas
\(\operatorname{rank}\mathcal F\ge n\ge3\). Hence the general fiber of \(\pi\) is \(\PP^k\).

Choose a general complete curve \(C\subset Y^\circ\), set
\(X_C:=X^\circ\times_{Y^\circ}C\), and suppress the pullback notation. We obtain a smooth \(\PP^k\)-bundle
\[
\pi:X_C\to C.
\]
By Tsen's theorem, there exists a vector bundle \(G\) on \(C\) such that
\[
X_C\simeq\PP_C(G).
\]
Moreover, \(\mathcal E|_{\PP^k}\simeq\OO_{\PP^k}(1)^{\oplus r}\) on every fiber. Set
\[
M:=\pi_*\left(\mathcal E\otimes\OO_{\PP(G)}(-1)\right).
\]
Cohomology and base change show that \(M\) is locally free of rank \(r\), and the natural evaluation morphism is an isomorphism
\[
\mathcal E\simeq\pi^*M\otimes\OO_{\PP(G)}(1).
\]
Consequently,
\[
\wedge^2\mathcal E
\simeq
\pi^*\wedge^2M\otimes\OO_{\PP(G)}(2).
\]
The inclusion
\(\wedge^2\mathcal E\subset\wedge^2T_{X_C/C}\) therefore gives a nonzero section of
\[
\wedge^2T_{X_C/C}
\otimes\pi^*\wedge^2M^*
\otimes\OO_{\PP(G)}(-2),
\]
contradicting Lemma \ref{vanish}.

Thus the assumption \(\rho(X)>1\) is impossible. Hence \(\rho(X)=1\), and the theorem follows from \cite{ross2010characterizations}. This proves Theorem \ref{rigid}, and therefore also Theorem \ref{maintheorem}.
\eproof

\section{Regular foliations whose exterior products are strictly nef}

In this section, we prove Theorem \ref{main} and Corollary \ref{coro}. The starting point is Ou's structure theorem \cite[Theorem 1.1]{ou2021foliations}: if a foliation $\mathcal F$ has nef $-K_{\mathcal F}$ and is regular (or has a compact leaf), then there exists a locally trivial fibration $f:X\to Y$ with rationally connected fibers and a foliation $\mathcal G$ on $Y$ such that $K_{\mathcal G}\equiv0$ and $\mathcal F=f^{-1}\mathcal G$. We also recall the definition of Shafarevich maps, which will be used to prove that the base of the fibration is of general type.

\bdefinition[\cite{kollar1995shafarevich}, Definition 3.5]
Let $X$ be a normal variety and let $H\subset \pi_1(X)$ be a normal subgroup. A normal variety $\operatorname{Sh}_H(X)$ together with a rational map
\[
\operatorname{sh}_X^H:X\dashrightarrow \operatorname{Sh}_H(X)
\]
is called the $H$-Shafarevich variety and the $H$-Shafarevich map of $X$ if the following conditions hold.

\noindent\textup{(1)} The map $\operatorname{sh}_X^H$ has connected fibers.

\noindent\textup{(2)} There are countably many proper closed subvarieties $D_i\subset X$ such that, for every closed subvariety $Z\subset X$ with $Z\not\subset \bigcup_i D_i$, the image $\operatorname{sh}_X^H(Z)$ is a point if and only if the image of $\pi_1(\tilde Z)\to \pi_1(X)$ has finite index in $H$, where $\tilde Z$ is the normalization of $Z$.
\edefinition

The existence and basic properties of Shafarevich maps are given by the following theorem.

\btheorem[\cite{kollar1995shafarevich}, Theorem 3.6]\label{shafarevich}
Let $X$ be a normal variety and let $H\subset \pi_1(X)$ be a normal subgroup. Then the $H$-Shafarevich map exists. Moreover, for every representative of the birational class of $\operatorname{Sh}_H(X)$, there are open subsets $X_0\subset X$ and $W_0\subset \operatorname{Sh}_H(X)$ such that $\operatorname{sh}_X^H$ is defined on $X_0$ and every fiber of $\operatorname{sh}_X^H|_{X_0}$ is closed in $X$. Furthermore, $\operatorname{sh}_X^H|_{X_0}$ is a topologically locally trivial fibration.
\etheorem

We will also use the following structure theorem for fibrations with nef relative anti-canonical divisor.

\btheorem[\cite{ccm2019}, \cite{matsumura2021structure}, Theorem 4.7]\label{locallyconstant}
Let $X$ be a normal variety with at worst klt singularities, and let $f:X\to Y$ be a fibration onto a smooth projective variety. Let $F$ be a general fiber of $f$. If $-K_{X/Y}$ is nef, then $f$ is a locally constant fibration induced by a representation $\pi_1(Y)\to \operatorname{Aut}(F)$. Moreover, there exists a very ample line bundle $A_F$ on $F$ which is $\pi_1(Y)$-linearized.
\etheorem

We recall two results from \cite{liu2020projective}, which will be used to prove hyperbolicity.

\blemma[\cite{liu2020projective}, Lemma 9.1]\label{degenerate}
Let $Z$ be a projective variety. If there exists a subgroup $G\subset \pi_1(Z)$ of finite index admitting a linear representation whose image is not virtually abelian, then every holomorphic map $\mathbb C\to Z$ is algebraically degenerate.
\elemma

\blemma[\cite{liu2020projective}, Lemma 9.2]\label{9.2}
Let $Z$ be a positive-dimensional projective variety, and let $\mathcal G$ be a flat vector bundle on $Z$ defined by a representation $\rho:\pi_1(Z)\to \operatorname{GL}_r(\mathbb C)$. If the image of $\rho$ is virtually solvable, then $\mathcal G$ is not strictly nef.
\elemma

We shall use the following elementary fact several times.

\blemma\label{extension}
Let
\[
0\longrightarrow E\longrightarrow F\longrightarrow G\longrightarrow 0
\]
be a short exact sequence of vector bundles on a projective variety $X$. Let $E'$ be a quotient vector bundle of $E$. Then there exists a quotient vector bundle $F'$ of $F$ fitting into a commutative diagram
\[
\begin{tikzcd}
0 \arrow[r] & E \arrow[r] \arrow[d] & F \arrow[r] \arrow[d] & G \arrow[r] \arrow[d,equal] & 0 \\
0 \arrow[r] & E' \arrow[r]          & F' \arrow[r]          & G \arrow[r]                 & 0 .
\end{tikzcd}
\]
\elemma

\bproof
The extension $0\to E\to F\to G\to 0$ is represented by an element $\tau\in \operatorname{Ext}^1(G,E)$. Let $E\to E'$ be the quotient map. The image of $\tau$ under the natural map $\operatorname{Ext}^1(G,E)\to \operatorname{Ext}^1(G,E')$ defines an extension
\[
0\longrightarrow E'\longrightarrow F'\longrightarrow G\longrightarrow 0.
\]
By functoriality of extensions, there is a morphism $F\to F'$ making the diagram commute. A diagram chase shows that $F\to F'$ is surjective.
\eproof

The next lemma is the criterion that will allow us to prove that the base of the fibration is canonically polarized and hyperbolic.

\blemma\label{hyper}
Let $Y$ be a smooth projective variety, and let $\rho:\pi_1(Y)\to \operatorname{GL}(k,\mathbb C)$ be a linear representation. Assume that for every positive-dimensional closed subvariety $Z\subset Y$, the image of the composition
\[
\pi_1(Z)\longrightarrow \pi_1(Y)\stackrel{\rho}{\longrightarrow}\operatorname{GL}(k,\mathbb C)
\]
is not virtually solvable. Then $Y$ is hyperbolic and $K_Y$ is ample.
\elemma

\bproof
We first prove hyperbolicity. Suppose that there exists a nonconstant holomorphic map $g:\mathbb C\to Y$. Let $Z\subset Y$ be the Zariski closure of $g(\mathbb C)$. By assumption, the image of $\pi_1(Z)\to \operatorname{GL}(k,\mathbb C)$ is not virtually solvable, hence in particular not virtually abelian. Lemma \ref{degenerate} implies that every holomorphic map $\mathbb C\to Z$ is algebraically degenerate, contradicting the definition of $Z$. Thus $Y$ is hyperbolic.

We next show that $Y$ is of general type. Let $\overline G$ be the Zariski closure of $\rho(\pi_1(Y))$, let $R$ be the solvable radical of $\overline G$, and let $H=\overline G/R$ be the semisimple quotient. Since $\rho(\pi_1(Y))$ is not virtually solvable, the induced representation $\rho_H:\pi_1(Y)\to H$ is nontrivial. After replacing $Y$ by a finite \'etale cover if necessary, we may assume that the image of $\rho_H$ is torsion-free.

Let $\operatorname{sh}_H:Y\dashrightarrow \operatorname{Sh}_H(Y)$ be the Shafarevich map associated with the kernel of $\rho_H$. By \cite[Th\'eor\`eme 1]{campana2015representations}, the Shafarevich variety $\operatorname{Sh}_H(Y)$ is of general type. We claim that $\operatorname{sh}_H$ is birational. If a positive-dimensional subvariety $Z\subset Y$ were contracted by $\operatorname{sh}_H$, then the image of $\pi_1(Z)\to H$ would be finite. Since the image of $\rho_H$ is torsion-free, this image would be trivial. Hence the image of $\pi_1(Z)$ under $\rho$ would be contained, up to finite index, in the inverse image of the solvable radical $R$, and therefore would be virtually solvable. This contradicts the assumption. Thus $\operatorname{sh}_H$ is birational, and $Y$ is of general type.

Since $Y$ is hyperbolic, it contains no rational curves. By the cone theorem, $K_Y$ is nef. Since $Y$ is of general type, $K_Y$ is nef and big. If $K_Y$ were not ample, the basepoint-free theorem would give a birational morphism from $Y$ to its canonical model with nonempty exceptional locus. Such an exceptional locus contains rational curves by the standard bend-and-break argument for birational contractions, contradicting hyperbolicity. Hence $K_Y$ is ample.
\eproof

We now show that the transverse part of the foliation must vanish.

\blemma\label{gzero}
Let $X$ be a projective manifold and let $\mathcal F$ be a regular foliation of rank $d$ on $X$. Suppose that there exists a locally constant fibration $f:X\to Y$ and a foliation $\mathcal G$ on $Y$ such that $c_1(\mathcal G)\equiv 0$ and $\mathcal F=f^{-1}\mathcal G$. If $\wedge^r\mathcal F$ is strictly nef for some $1\le r\le d$, then $\mathcal G=0$ and $\mathcal F=T_{X/Y}$.
\elemma

\bproof
Assume by contradiction that $\mathcal G\neq 0$. Since $\mathcal F=f^{-1}\mathcal G$, we have an exact sequence
\[
0\longrightarrow T_{X/Y}\longrightarrow \mathcal F\longrightarrow f^*\mathcal G\longrightarrow 0.
\]
Taking exterior powers gives a quotient $\wedge^r\mathcal F\to \wedge^r f^*\mathcal G$. Since $\wedge^r\mathcal F$ is strictly nef, every nonzero quotient of it is strictly nef. If $f$ is not the identity and $\wedge^r\mathcal G\neq 0$, then the pullback $f^*\wedge^r\mathcal G$ has degree zero on every curve contained in a fiber of $f$, and hence cannot be strictly nef. Thus either $f$ is the identity, or $\wedge^r\mathcal G=0$.

First suppose that $f$ is the identity. Then $\mathcal F=\mathcal G$, and $\wedge^r\mathcal G$ is strictly nef. Since $c_1(\mathcal G)\equiv 0$, the bundle $\wedge^r\mathcal G$ is numerically flat. By \cite{demailly1994compact}, Theorem on page 3, it admits a filtration whose graded quotients are Hermitian flat vector bundles. Hence $\wedge^r\mathcal G$ has a nonzero flat quotient bundle $\mathcal Q$, which is strictly nef.

Let $\rho_{\mathcal Q}:\pi_1(Y)\to \operatorname{GL}(q,\mathbb C)$ be the representation associated with $\mathcal Q$. For every positive-dimensional closed subvariety $Z\subset Y$, the restriction $\mathcal Q|_Z$ is still strictly nef and flat. By Lemma \ref{9.2}, the image of $\pi_1(Z)$ under $\rho_{\mathcal Q}$ is not virtually solvable. Hence Lemma \ref{hyper} applies and shows that $K_Y$ is ample.

By the Aubin--Yau theorem, the tangent bundle $T_Y$ is semistable with respect to $K_Y$. Since $\mathcal G\subset T_Y$ and $c_1(\mathcal G)\equiv 0$, we have
\[
\mu_{K_Y}(\mathcal G)=0>\mu_{K_Y}(T_Y),
\]
which contradicts semistability. Therefore the case $f=\operatorname{id}_Y$ is impossible.

We may now assume that $f$ is not the identity. Then the preceding observation gives $\wedge^r\mathcal G=0$, hence $\operatorname{rank}\mathcal G<r$. Since $c_1(\mathcal G)\equiv 0$, we have
\[
-K_{\mathcal F}\equiv -K_{X/Y}.
\]
Moreover,
\[
\det(\wedge^r\mathcal F)\simeq (-K_{\mathcal F})^{\otimes \binom{d-1}{r-1}}
\]
is nef, because $\wedge^r\mathcal F$ is strictly nef. Hence $-K_{X/Y}$ is nef. By Theorem \ref{locallyconstant}, we may choose a very ample $\pi_1(Y)$-linearized line bundle $A_F$ on a fiber $F$ of $f$. Let
\[
\rho:\pi_1(Y)\to \operatorname{GL}(V)
\]
be the induced representation, where $V=H^0(F,A_F)$.

We claim that for every positive-dimensional closed subvariety $Z\subset Y$, the image of $\pi_1(Z)\to \pi_1(Y)\stackrel{\rho}{\to}\operatorname{GL}(V)$ is not virtually solvable. Suppose otherwise. After replacing $Z$ by a finite \'etale cover, Borel's fixed point theorem gives a fixed point $p\in F$ for the induced action on $F$. This fixed point produces a section $W\subset X_Z:=X\times_Y Z$ of the base-changed fibration $X_Z\to Z$.

The restriction $T_{X_Z/Z}|_W$ is a flat vector bundle associated with the differential of the monodromy action at $p$. Since the monodromy image is virtually solvable, after another finite \'etale cover of $Z$ the bundle $T_{X_Z/Z}|_W$ has a filtration by flat line bundles. Because $\operatorname{rank}\mathcal G<r$, we can choose a flat quotient bundle $Q$ of $T_{X_Z/Z}|_W$ of rank $r-\operatorname{rank}\mathcal G$.

Restricting the exact sequence defining $\mathcal F$ to $W$, and applying Lemma \ref{extension}, we obtain a quotient bundle $Q'$ of $\mathcal F|_W$ fitting into an exact sequence
\[
0\longrightarrow Q\longrightarrow Q'\longrightarrow \mathcal G|_Z\longrightarrow 0.
\]
Thus $\operatorname{rank}Q'=r$, and
\[
\det Q'\equiv \det Q+\det \mathcal G|_Z\equiv 0.
\]
The quotient map $\mathcal F|_W\to Q'$ induces a quotient $\wedge^r\mathcal F|_W\to \det Q'$. Since $\wedge^r\mathcal F$ is strictly nef, the line bundle $\det Q'$ must be strictly nef, contradicting $\det Q'\equiv 0$. This proves the claim.

By the claim and Lemma \ref{hyper}, the base $Y$ is hyperbolic and $K_Y$ is ample. As before, the Aubin--Yau theorem implies that $T_Y$ is semistable with respect to $K_Y$. Since $\mathcal G\subset T_Y$ and $c_1(\mathcal G)\equiv 0$, this contradicts semistability. Hence $\mathcal G=0$, and therefore $\mathcal F=T_{X/Y}$.
\eproof

We can now prove the main theorem.

\btheorem\label{strictlyneffo}
Let $\mathcal F$ be a regular foliation of rank $d$ on a projective manifold $X$. If $\wedge^r\mathcal F$ is strictly nef for some $1\le r<d$, then there exists a locally constant fibration
\[
f:X\to Y
\]
such that the fibers are Fano manifolds and the base $Y$ is a canonically polarized hyperbolic projective manifold. Moreover, the foliation $\mathcal F$ is induced by the fibration $f$. In particular, $\mathcal F$ is algebraically integrable.
\etheorem

\bproof
By Ou's structure theorem, there exists a locally trivial fibration $f:X\to Y$ with rationally connected fibers and a foliation $\mathcal G$ on $Y$ with $K_{\mathcal G}\equiv 0$ such that $\mathcal F=f^{-1}\mathcal G$. Since
\[
-K_{\mathcal F}\equiv -K_{X/Y}
\]
and
\[
\det(\wedge^r\mathcal F)\simeq (-K_{\mathcal F})^{\otimes \binom{d-1}{r-1}},
\]
the strict nefness of $\wedge^r\mathcal F$ implies that $-K_{X/Y}$ is nef. By Theorem \ref{locallyconstant}, the fibration $f$ is locally constant and admits a $\pi_1(Y)$-linearized very ample line bundle on the fiber.

By Lemma \ref{gzero}, we have $\mathcal G=0$. Hence $\mathcal F=T_{X/Y}$. Therefore $\wedge^rT_{X/Y}$ is strictly nef, and its restriction to every fiber $F$ of $f$ gives that $\wedge^rT_F$ is strictly nef. Since $1\le r<d=\dim F$, Watanabe's result \cite{watanabe2022positivity}, Remark after Proposition 1.4, implies that every fiber $F$ is Fano.

It remains to prove that $Y$ is hyperbolic and that $K_Y$ is ample. Let $A_F$ be a very ample $\pi_1(Y)$-linearized line bundle on a fiber $F$, and let
\[
\rho:\pi_1(Y)\to \operatorname{GL}(H^0(F,A_F))
\]
be the induced representation. We show that the hypothesis of Lemma \ref{hyper} is satisfied.

Let $Z\subset Y$ be a positive-dimensional closed subvariety, and consider the base change $X_Z:=X\times_Y Z\to Z$. If the image of $\pi_1(Z)$ under $\rho$ were virtually solvable, then, after replacing $Z$ by a finite \'etale cover, Borel's fixed point theorem would give a fixed point $p\in F$. This fixed point would induce a section $W\subset X_Z$. As in the proof of Lemma \ref{gzero}, the bundle $T_{X_Z/Z}|_W$ would be flat with virtually solvable monodromy, and hence, after a finite \'etale cover, would admit a filtration by flat line bundles. We could then choose a flat quotient bundle $Q$ of $T_{X_Z/Z}|_W$ of rank $r$.

The quotient $T_{X_Z/Z}|_W\to Q$ would induce a quotient
\[
\wedge^rT_{X_Z/Z}|_W\to \det Q.
\]
But $\det Q\equiv 0$, whereas $\wedge^rT_{X_Z/Z}|_W$ is strictly nef. This is impossible. Therefore the image of $\pi_1(Z)$ under $\rho$ is not virtually solvable for every positive-dimensional closed subvariety $Z\subset Y$.

By Lemma \ref{hyper}, the base $Y$ is hyperbolic and $K_Y$ is ample. This completes the proof.
\eproof
\bcorollary\label{coro}
Let $\mathcal F$ be a regular foliation of rank $d\ge 3$ on a projective manifold $X$. If $\wedge^2\mathcal F$ is strictly nef, then there exists a locally constant fibration
\[
f:X\to Y
\]
whose fibers are projective spaces or smooth hyperquadrics, where $Y$ is a canonically polarized hyperbolic projective manifold. Moreover, $\mathcal F$ is the foliation induced by the fibration. In particular, $\mathcal F$ is algebraically integrable.
\ecorollary
\bproof[Proof of Corollary \ref{coro}]
By Theorem \ref{strictlyneffo}, the foliation $\mathcal F$ is induced by a locally constant fibration $f:X\to Y$, and $\mathcal F=T_{X/Y}$. Hence for every fiber $F$ of $f$, the bundle $\wedge^2T_F$ is strictly nef. By \cite{li2019projective}, Theorem 1.5, the fiber $F$ is isomorphic either to a projective space or to a smooth hyperquadric. The remaining assertions follow directly from Theorem \ref{strictlyneffo}.
\eproof
\bibliographystyle{plain}
\bibliography{references}

@article{andreatta2001manifolds,
  title={On manifolds whose tangent bundle contains an ample subbundle},
  author={Andreatta, Marco and Wi{\'s}niewski, Jaros{\l}aw A},
  journal={Inventiones Mathematicae},
  volume={146},
  number={1},
  pages={209--217},
  year={2001},
  publisher={Springer Berlin Heidelberg Berlin/Heidelberg}
}

@article{mori1979projective,
  title={Projective manifolds with ample tangent bundles},
  author={Mori, Shigefumi},
  journal={Annals of Mathematics},
  volume={110},
  number={3},
  pages={593--606},
  year={1979},
  publisher={JSTOR}
}

@article{siu1980compact,
  title={Compact {K}{\"a}hler manifolds of positive bisectional curvature},
  author={Siu, Yum-Tong and Yau, Shing-Tung},
  journal={Inventiones Mathematicae},
  volume={59},
  number={2},
  pages={189--204},
  year={1980},
  publisher={Springer-Verlag Berlin/Heidelberg}
}

@article{demailly1994compact,
  title={Compact complex manifolds with numerically effective tangent bundles},
  author={Demailly, Jean-Pierre and Peternell, Thomas and Schneider, Michael},
  journal={Journal of Algebraic Geometry},
  volume={3},
  number={2},
  pages={295--346},
  year={1994},
  publisher={Providence, RI: University Press, c1992-}
}

@article{druel2013characterizations,
  title={Characterizations of projective spaces and hyperquadrics},
  author={Druel, St{\'e}phane and Paris, Matthieu},
  journal={Asian Journal of Mathematics},
  volume={17},
  number={4},
  pages={583--596},
  year={2013},
  publisher={International Press of Boston}
}

@article{liu2019characterization,
  title={Characterization of projective spaces and-bundles as ample divisors},
  author={Liu, Jie},
  journal={Nagoya Mathematical Journal},
  volume={233},
  pages={155--169},
  year={2019},
  publisher={Cambridge University Press}
}

@article{mok1988uniformization,
  title={The uniformization theorem for compact {K}{\"a}hler manifolds of nonnegative holomorphic bisectional curvature},
  author={Mok, Ngaiming},
  journal={Journal of differential geometry},
  volume={27},
  number={2},
  pages={179--214},
  year={1988},
  publisher={Lehigh University}
}

@article{ross2010characterizations,
  title={Characterizations of projective spaces and hyperquadrics via positivity properties of the tangent bundle},
  author={Ross, Kiana},
  journal={arXiv preprint arXiv:1012.2043},
  year={2010}
}

@article{araujo2007cohomological,
  title={Cohomological characterizations of projective spaces and hyperquadrics},
  author={Araujo, Carolina and Druel, St{\'e}phane and Kov{\'a}cs, S{\'a}ndor J},
  journal={arXiv preprint arXiv:0707.4310},
  year={2007}
}

@article{liu2020projective,
  title={Projective manifolds whose tangent bundle contains a strictly nef subsheaf},
  author={Liu, Jie and Ou, Wenhao and Yang, Xiaokui},
    journal = {J. Algebraic Geom.},
    volume = {33},
      year = {2024},
     pages = {1--53},
}

@article{campana1991projective,
  title={Projective manifolds whose tangent bundles are numerically effective},
  author={Campana, Fr{\'e}d{\'e}ric and Peternell, Thomas},
  journal={Mathematische Annalen},
  volume={289},
  pages={169--187},
  year={1991},
  publisher={Springer}
}

@article{yasutake2012second,
  title={On the second exterior power of tangent bundles of {F}ano fourfolds with {P}icard number $ \rho(X) \geq 2$},
  author={Yasutake, Kazunori},
  journal={arXiv preprint arXiv:1212.0685},
  year={2012}
}

@article{cho1995smooth,
  title={Smooth projective varieties with the ample vector bundle $\wedge^2{T}_{X}$ in any characteristic},
  author={Cho, Koji and Sato, Eiichi},
  journal={Journal of Mathematics of Kyoto University},
  volume={35},
  number={1},
  pages={1--33},
  year={1995},
  publisher={Duke University Press}
}

@article{li2019projective,
  title={On projective varieties with strictly nef tangent bundles},
  author={Li, Duo and Ou, Wenhao and Yang, Xiaokui},
  journal={Journal de Math{\'e}matiques Pures et Appliqu{\'e}es},
  volume={128},
  pages={140--151},
  year={2019},
  publisher={Elsevier}
}

@article{watanabe2021positivity,
  title={Positivity of the second exterior power of the tangent bundles},
  author={Watanabe, Kiwamu},
  journal={Advances in Mathematics},
  volume={385},
  pages={107757},
  year={2021},
  publisher={Elsevier}
}

@article{watanabe2022positivity,
  title={Positivity of the exterior power of the tangent bundles},
  author={Watanabe, Kiwamu},
  journal={Proc. Japan Acad. Ser. A Math. Sci. },
  year={2023},
volume={99(10)},
pages={77-80}
}

@article{cao2017decomposition,
  title={A decomposition theorem for projective manifolds with nef anticanonical bundle},
  author={Cao, Junyan and H{\"o}ring, Andreas},
  journal={arXiv preprint arXiv:1706.08814},
  year={2017}
}

@book{lazarsfeld2017positivity,
  title={Positivity in algebraic geometry I: Classical setting: line bundles and linear series},
  author={Lazarsfeld, Robert K},
  volume={48},
  year={2017},
  publisher={Springer}
}

@article{liu2021algebraic,
  title={Algebraic fibre spaces with strictly nef relative anti-log canonical divisor},
  author={Liu, Jie and Ou, Wenhao and Wang, Juanyong and Yang, Xiaokui and Zhong, Guolei},
  journal={arXiv preprint arXiv:2111.05234},
  year={2021}
}

@article{greb2016movable,
  title={Movable curves and semistable sheaves},
  author={Greb, Daniel and Kebekus, Stefan and Peternell, Thomas},
  journal={International Mathematics Research Notices},
  volume={2016},
  number={2},
  pages={536--570},
  year={2016},
  publisher={Oxford University Press}
}

@article{ou2021foliations,
  title={Foliations whose first class is nef},
  author={Ou, Wenhao},
  journal={to appear on Journal of the European Mathematical Society},
  year={2021}
}

@article{ccm2019,
  title={Projective klt pairs with nef anti-canonical divisor},
  author={Campana, Fr{\'e}d{\'e}ric and Cao, Junyan and Matsumura, Shinichi},
  journal={Algebraic Geometry},
  volume={8},
  pages={430--464},
  year={2021}
}

@article{matsumura2021structure,
  title={Structure theorem for projective klt pairs with nef anti-canonical divisor},
  author={Matsumura, Shinichi and Wang, Juanyong},
  journal={arXiv preprint arXiv:2105.14308},
  year={2021}
}

@article{campana2015representations,
  title={Repr{\'e}sentations lin{\'e}aires des groupes k{\"a}hl{\'e}riens: factorisations et conjecture de Shafarevich lin{\'e}aire},
  author={Campana, Fr{\'e}deric and Claudon, Beno{\^\i}t and Eyssidieux, Philippe},
  journal={Compositio Mathematica},
  volume={151},
  number={2},
  pages={351--376},
  year={2015},
  publisher={London Mathematical Society}
}

@book{kollar1995shafarevich,
  title={Shafarevich maps and automorphic forms},
  author={Koll{\'a}r, J{\'a}nos},
  year={1995},
  publisher={Princeton University Press}
}
\end{document}